\documentclass [11pt] {article}
\usepackage{amsfonts}
\usepackage{amssymb}
\newcommand{\qed} {\hspace {0.1in} \rule {1.5mm} {3.5mm}}

\newtheorem{lemma}{Lemma}[section]
\newtheorem{corollary}{Corollary}[section]

\newtheorem{theorem}{Theorem}

\newtheorem{proposition}{Proposition}[section]

\def\limn{\lim_{n\to\infty}}
\def\gseq{\{G_n\}_{n=1}^\infty}

\def\e{\epsilon}

\def\limo{\lim_{\omega}}

\def\dim{{\rm dim}}

\def\<{\langle}
\def\>{\rangle}

\def\proof{\smallskip\noindent{\it Proof.} }

\def\bZ{{\mathbb Z}}
\def\bG{{\bf G}}
\def\bA{{\bf A}}
\def\bX{{\bf X}}
\def\bB{{\bf B}}
\def\bE{{\bf E}}
\def\bH{{\bf H}}

\def\bK{{\bf K}}
\def\bbbN{{\bf N}}

\def\bQ{{\bf Q}}
\def\bO{{\bf O}}
\def\bbN{{\bf N}}
\def\bR{{\mathbb R}}
\def\bN{{\mathbb N}}

\def\bS{{\bf S}}

\def\cB{\mbox{$\cal B$}}

\def\cH{\mbox{$\cal H$}}
\def\cG{\mbox{$\cal G$}}

\def\cF{\mbox{$\cal F$}}

\def\cB{\mbox{$\cal B$}}

\def\cP{\mbox{$\cal P$}}
\def\cR{\mbox{$\cal R$}}

\def\cM{\mbox{$\cal M$}}

\def\Fo{F$\mbox{\o}$lner}

\def\gseq{\{G_n\}^\infty_{n=1}}
\def\gd{{\bf Gr}_d}
\def\grd{\mbox{Graph}_d }
\def\bgd{{\bf BGr}_d}
\def\bgf{{\bf BGr}_f}
\def\grf{\mbox{Graph}^f_d }
\def\f2{{\bf F}^\infty_2}
\def\hmu{\hat{\mu}_{\bG}}
\def\tmu{\widetilde{\mu}_{\bG}}
\def\uz{\underline{z}}
\def\bp{{\bf p}}
\def\bq{{\bf q}}
\def\bmu{\mu_{\bG}}
\def\xg{\bX_{\bG}}

\begin{document}
\title{Parameter testing in bounded degree graphs of subexponential growth}
\author{{\sc G\'abor Elek}
\footnote {The Alfred Renyi Mathematical Institute of
the Hungarian Academy of Sciences, P.O. Box 127, H-1364 Budapest, Hungary.
email:elek@renyi.hu, Supported by OTKA Grants T 049841 and T 037846}}
\date{}
\maketitle \vskip 0.2in \noindent
\vskip 0.2in \noindent{\bf Abstract.} Parameter testing algorithms are using
constant number of queries to estimate the value of a certain parameter of a
very large finite graph. It is well-known that graph parameters such as the
independence ratio or the edit-distance from $3$-colorability are not testable
in bounded degree graphs. We prove, however, that these and several other
interesting graph parameters are testable in bounded degree graphs of
subexponential growth. 

\vskip 0.3in

\newpage
\tableofcontents
\newpage
\section{Introduction}
\subsection{Dense graph sequences}
The main motivation for our paper is to develop a theory analogous to that
recently developed
for dense graph sequences \cite{BCh1},\cite{BCh2},\cite{LSZ}.
First let us recall some basic notions.
A sequence of finite simple graphs $\bG=\gseq$, $|V(G_n)|\to\infty$
 is called {\bf convergent}
if for any finite simple graph $F$, $\limn t(F,G_n)$ exists where
$$t(F,G)=\frac{|\mbox{hom}\,(F,G)|}{|V(G)|^{|V(F)|}}$$
is the probability that a random map from $V(F)$ into $V(G)$ is a graph
homomorphism. The convergence structure above defines a metrizable
compactification of the sets of finite graphs. The {\bf limit objects}
of the graph sequences were first introduced in \cite{LSZ}. They are
measurable symmetric functions
$$W:[0,1]\times[0,1]\to [0,1]\,.$$
A graph sequence $\gseq$ converges to $W$
if for every finite simple graph $F$,
$$ \limn t(F,G_n)=\int_{[0,1]^{V(F)}}\prod_{(i,j)\in E(F)} W(x_i,x_j)
dx_1 dx_2\dots dx_{|V(F)|}\,.$$
For any such function $W$ one can find a graph sequence $\gseq$
converging to $W$ and conversely for any graph sequence $\gseq$ there
exists a measurable function $W$ such that the sequence converges to $W$.
Consequently, the boundary points of the compactification can be identified
with equivalence classes of such measurable functions \cite{LSZ}. Note
that if $\gseq$ is a {\bf sparse} sequence with
$\limn \frac{|E(G_n)|}{|V(G_n)|^2}=0$, then $\gseq$ in fact converges to
the zero function.

\noindent
A {\bf graph parameter} is a real function on the sets of finite
simple graphs that is invariant under graph isomorphims. A parameter
$\phi$ is {\bf continuous} 
 if $\limn \phi(G_n)$ exists for any convergent sequence $\gseq$.
It was shown by Fischer and Newman \cite{FN}
 that continuous graph parameters are
exactly the ones that are {\bf testable} by random samplings.
It has been proved first in \cite{AS}
and then later in \cite{LSZ} that the
edit-distance from a hereditary graph property is a continuous graph
parameter.

\subsection{Bounded degree graphs}
Let $d\geq 2$ be a positive integer and let $\grd$ be the set of finite
graphs $G$ (up to isomorphisms) such that
$\mbox{deg}(x)\leq d$ for any $x\in V(G)$. The notion of {\bf weak
  convergence} for the class $\grd$ was introduced by Benjamini and
Schramm \cite{BS}. Let us start with some definitions.
A rooted $(r,d)$-ball is a finite, simple, connected graph $H$ such that
\begin{itemize}
\item $\mbox{deg}(y)\leq d$ if $y\in V(H)$\,.
\item $H$ has a distinguished vertex $x$ (the root).
\item $d_G(x,y)\leq r$ for any $y\in V(H)$.
\end{itemize}
For $r\geq 1$, we denote by $U^{r,d}$ the finite set of rooted isomorphism
classes of rooted $(r,d)$-balls. Let $G(V,E)$ be a finite graph with vertex
degree bound $d$. For $\alpha\in U^{r,d}$, $T(G,\alpha)$ denotes the
set of vertices $x\in V(G)$ such that there exists a rooted isomorphism
between $\alpha$ and the rooted $r$-ball $B_r(x)$ around $x$.
 Set $p_G(\alpha):=\frac{|T(G,\alpha)|}
{|V(G)|}\,.$ Thus we associated to $G$ a probability distribution on $U^{r,d}$
for any $r\geq 1$. Let $\bG=\gseq\subset \grd$ be a sequence of finite
simple graphs such that $\limn |V(G_n)|=\infty$. Then $\bG$ is called
{\bf weakly convergent} if for any $r\geq 1$ and $\alpha\in U^{r,d}$,
$\limn p_{G_n}(\alpha)$ exists. The convergence structure above
defines a metrizable compactification of $\grd$ in the following way.
Let $\alpha_1,\alpha_2,\dots$ be an enumeration of the elements of
$\cup^\infty_{r=1} U^{r,d}$.
For a graph $G$ we associate a sequence
$$s(G)=\{\frac{1}{|V(G)|}, p_G(\alpha_1), p_G(\alpha_2),\dots\}\in
[0,1]^{\bN}\,.$$
 By definition,
$\{G_n\}^\infty_{n=1}$ is weakly convergent if and only if 
$\{s(G_n)\}^\infty_{n=1}$ converge pointwise. We consider the closure of
$s(\grd)$ in the compact space $[0,1]^{\bN}$. This set can be viewed as the
compactification of $\grd$.
Again, a graph parameter
$\phi:\grd\to\bR$ is called continuous if $\limn \phi(G_n)$ exists for any
weakly convergent sequence. Equivalently, $\phi$ is continuous if it extends
continuously to the compactification above. 

\subsection{Hyperfinite graph classes}
Hyperfinite graph classes were introduced in \cite{Elek2} and studied in depth
in \cite{Sch},\cite{BSS}. Also, under the name of non-expanding bounded degree
graph classes they were studied in \cite{CSS} as well.
A class $\mathcal{H}\subset \grd$ is called {\bf hyperfinite} if for
any $\e>0$ there exists $K>0$ such that if $G\in\mathcal{H}$ then one can
delete $\e|E(G)|$ edges from $G$ in such a way that all the components
in the remaining graph $G'$ have size at most $K$. Planar graphs, graphs with
bounded treewidth or, in general, all the minor-closed graph classes are
hyperfinite \cite{BSS}. Let $f:\bN\to\bN$ be a function of subexponential
growth. That is, for any $\delta>0$ there exists $C_\delta>0$ such that
for all $n\geq 1$: $ f(n)\leq C_\delta(1+\delta)^n$. The class $\grf$ consists
of graphs $G\in\grd$ such that $|B_r(x)|\leq f(r)$ for each $x\in V(G)$.
The classes $\grf$ are also hyperfinite \cite{Elek2}.
We will call a graph parameter $\phi$ continuous on $\grf$ if 
$\lim_{n\to\infty} \phi(G_n)$ exists whenever $\{G_n\}^\infty_{n=1}\subset
\grf$ is a weakly convergent sequence. Equivalently, $\phi$ is continuous on
$\grf$ if it extends continuously to the closure of $\grf$ in $[0,1]^\bN$.

\subsection{Union-closed monotone properties}
Let $\cP\subset\grd$. We say that $\cP$ is union-closed monotone graph class
(or being in $\cP$ is
a union-closed monotone property) 
if the following conditions are satisfied:
\begin{itemize}
\item if $|E(A)|=0$ then $A\in \cP$
\item if $A\in \cP$ and $B\subset A$ is a subgraph, then $B\in \cP$
(we consider spanning subgraphs, that is if $B\subset A$ then
$V(B)=V(A)$)
\item if $A\in \cP$ and $B\in \cP$ then the disjoint union of $A$ and $B$
is also in $\cP$. 
\end{itemize}
Let us list some union-closed monotone graph classes :
\begin{itemize}
\item planar graphs
\item bipartite graphs
\item $k$-colorable graphs
\item graphs that are not containing some fixed graph $H$.
\end{itemize}
If $G$ and $H$ are finite graphs with the same vertex set $V$ then
their {\bf edge-distance} is defined as
$$d_e(G,H):=\frac{|E(G)\triangle E(H)|}{|V|}\,.$$
The {\bf edit-distance} from a class $\cP$ is defined as
$$d_e(G,\cP)=\inf_{V(H)=V(G),\, H\in\cP} d_e(G,H)\,.$$
It is important to note that $d_e(*,\cP)$ is not continuous on
$\grd$ even for such a simple class as the set of bipartite graphs $Bip$.
Indeed, Bollob\'as \cite{Bol} constructed a large girth sequence of
cubic graphs such that $d_e(G_n,Bip)>\e>0$ for any $n\geq 1$. On the
other hand there are bipartite large girth sequences of cubic graphs.
Since by the definition of weak convergence all sequences of cubic graphs
with large girth converge to the same elements of the compactification of
$\grd$, $d_e(*,\cP)$ is not continuous.
 We shall see, however, that in
the class $\grf$ the graph parameter $d_e(*,\cP)$ is continuous if $\cP$
is a union-closed monotone property (Theorem \ref{theo2}). 
We also prove that continuous graph
parameters are effectively testable via random samplings (Theorem \ref{test}).

\subsection{Continuous graph parameters in $\grf$} 
In Section \ref{Ind} we prove that the independence ratio as well as
the matching ratio are continuous parameters for the class $\grf$ (Theorem
\ref{utolsotetel}). We
also prove that the log-partition functions associated to independent subsets
 resp.
to matchings are continuous graph parameters in $\grf$. This shows that for
certain aperiodic graphs such as the Penrose tilings, in which all
neighbourhood patterns can be seen in a given frequency, the thermodynamical
limit of the log-partition functions exists. Such results are
well-known for lattices. We also show a similar convergence result
for the integrated density of states for discrete Schr\"odinger
operators with random potentials extending some recent results in \cite{LS}
and \cite{LV}(Theorem \ref{int}).

\subsection{The main theorem}
It is known \cite{AK} that there exists $\delta>0$ such that to construct
an independent set that
approximates the size of a maximum independent set within an
error of $\delta |V(G)|$
in a $3$-regular graph $G$ is NP-hard. The situation is dramatically
different in the case of planar graphs. For any fixed $\delta>0$ 
there exists a polynomial time algorithm to construct
an independent set that approximates the size of the maximum
independent set within an error of $\delta |V(G)|$ for cubic
planar graphs $G$
\cite{Baker} (note that finding a maximum independent
set in a planar cubic graph is still NP-hard). First, using a
polynomial time algorithm
one can delete $\frac {\delta |E(G)|}{3}$ edges from $G$ to obtain a graph
$G'$ with components of size at most $K(\delta)$. For each component of $G'$
one can find a maximum independent set in $L(\delta)$ steps. Obviously, the
union
of these sets can not be smaller in size that the maximum independent set in
$G$. If we delete all the vertices from the union that are on some previously
deleted edges, then we get an independent subset of the original graph $G$.
Since the number of deleted vertices is at most $\delta |V(G)|$  we obtained
an approximation of the maximum independent set within an error of $\delta
|V(G)|$. 

\vskip0.1in
\noindent
How can we use this idea for constant-time algorithms ?
Let $f$ be a function of subexponential growth and $G\in\grf$. Fix $\e>0$.
Since $\grf$ is a hyperfinite class one can delete $\epsilon |E(G)|$
edges from $G$ to obtain a graph $G'$ with components of size at most
$K(\epsilon)$. Let $A(d,K(\e))$ be the finite
set of all finite connected graphs
of size at most $K(\epsilon)$. If for each $H\in A(d,K(\e))$ someone tells
us how many components of $G'$ are isomorphic to $H$ we can calculate
the size of the maximum independent set in $G'$. What we need is to test
the following value : the number of components in $G'$ isomorphic to $H$
divided by
$|V(G)|$. Unfortunately, this is not a well-defined graph parameter since 
there are
many ways to delete edges from $G$ to obtain graphs with small components.
Informally speaking, what we need to show is that if 
two graphs $G_1, G_2\in\grf$ are close to each other
in terms of local neighborhood statistics, then one can delete edges from $G_1$
resp. $G_2$ in such a way that in the remaining graphs $G'_1$ resp. $G'_2$ the
ratios of $H$-components are close to each other for any fixed $H\in
A(d,K(\e))$. That is exactly what we prove in our main theorem, which is the 
main tool of our paper.
\begin{theorem}\label{main}
Let $\bG=\gseq\subset\grf$ be
a weakly convergent sequence of finite graphs.
Then for any $\e>0$ there exists a constant $K>0$ and also, for all
connected simple graphs $H\in\grf$ with $|V(H)|\leq K$ a real constant
$c_H$ such that for any $n\geq 1$ one can remove $\e|E(G_n)|$ edges from $G_n$
satisfying the following conditions:
\begin{itemize}
\item The number of vertices in each component of the remaining
graph $G'_n$ is not greater than $K$.
\item If $V^n_H\subseteq V(G_n)$ is the
set of vertices that are contained in a component of $G_n'$ isomorphic
to $H$ then
$$\limn\frac{|V_H^n|}{|V(G_n)|}=c_H\,.$$
\end{itemize}
\end{theorem}
Note that the second condition is equivalent to saying that $\{G_n'\}$
is a convergent sequence.
In order to prove the theorem we combine the limit object method
of Benjamini and Schramm \cite{BS} and the non-standard 
analytic technique
developed in \cite{ESZ}.

\section{The canonical limit object}
\subsection{Hyperfinite graphings} \label{hyperfinite}
In this subsection we briefly recall the basic properties of graphings
(graphed equivalence relations) \cite{Kech}.
 Let $\f2$ be the free
product of countably many copies of the cycle group of order two.
Thus 
$$\f2=\langle \{s_i\}_{i\in \bN} | s^2_i=1\rangle $$
is a presentation of the group $\f2$, where the
$s_i$'s are generators of order two. Suppose that the edges of a simple
graph $H$ (finite or infinite) are coloured by natural numbers
{\bf properly}, that is, any two edges having a common vertex are coloured
differently.
Then, the colouring induces an action of $\f2$ on the vertex set
$V(H)$ in the following way:
\begin{itemize}
\item $s_i(x)=y$ if $e=(x,y)\in E(H)$ and $e$ is
coloured by $i$.
\item $s_i(x)=x$ if no edge incident to $x$ is coloured by $i$.
\end{itemize}
\noindent
We regard graphings as the measure
theoretical analogues of $\bN$-coloured graphs.
Let $(X,\mu)$ be a probability measure space with
a measure-preserving action of $\f2$ that is not necessarily free such that
if $s_i(p)=q\neq p$ and $s_j(p)=q$ then $i= j$. Let $E\subset X\times X$
be the set of pairs $(p,q)$ such that $\gamma(p)=q$ for some $\gamma\in \f2$.
Thus $E$ is the {\bf measurable equivalence relation} induced by the
$\f2$-action. Connect the points $p\in X, q\in X, p\neq q$
 by an edge of colour $i$ if
$s_i(p)=q$ for some generator element $s_i$. Thus we obtain a 
properly $\bN$-coloured graph
with a measurable structure, the {\bf graphing} $\cG$. If $p\in X$
then $\cG_p$ denotes the component of $\cG$ containing $p$.
In this paper we consider only bounded degree graphings, that is, 
graphings for which all the degrees of the vertices are bounded by
a certain constant $d$. Thus, for any $p\in X$
the number of generators $\{s_i\}$ which do not fix $p$ is at most $d$.
The edge-set of the graphing $\cG$, $E(\cG)$ has a natural measure
space structure as well.
Let $i\geq 1$ and $A\subseteq X$ be a measurable subset of vertices such 
that
\begin{itemize} \item
If $a\in A$ then $s_i(a)\in A$, $s_i(a)\neq a$.
\end{itemize}
Let $\hat{A_i}\subseteq E(\cG)$ be the set of edges such that
their endpoints belong to $A$. Then we call $\hat{A_i}$
a measurable edge-set of colour $i$.  These measurable edge-sets
form a $\sigma$-algebra with a measure $\mu_{E_i}$,
$$\mu_{E_i}(\hat{A_i})=\frac{1}{2}\mu(A)\,.$$
Clearly, the $\sigma$-algebra above contains the set $E_i$ consisting of all
edges coloured by $i$.
Then $E(\cG)=\cup^\infty_{i=1} E_i$. The set
$M\subset E(\cG)$ is measurable if for all $i$ $M\cap E_i$ is measurable
and
$$\mu_E(M)=\sum^\infty_{i=1}\mu_{E_i} (M\cap E_i)\,.$$
A {\bf measurable subgraphing} $\cH\subseteq\cG$ is
a measurable subset of $E(\cG)$ such that the components of $\cH$
are induced subgraphs of $\cG$. A subgraphing $\cH$ is called component-finite
if all of its components are finite graphs. It is easy to see that if
$\cH\subset \cG$ is a component-finite subgraphing and 
$F$ is a finite connected
simple graph, then
$$\cH_F=\{p\in X\,\mid\, \cH_p\cong F\}$$
is measurable and the span of $\cH_F$ is a component-finite subgraphing
having components isomorphic to $F$.

\noindent
The graphing $\cG$ is called {\bf hyperfinite} if there exist
component-finite subgraphings $\cH_1\subset\cH_2\subset\dots$ such that
$$\limn \mu_E(E(\cG)\backslash E(\cH_n))=0\,.$$
Suppose that $f:\bN\to \bR$ is a function of subexponential
growth and $|B_r(x)|\leq f(r)$ for all the balls of radius $r$ in
$\cG$. Then we call $\cG$ a graphing of subexponential growth.
By the result of Adams and Lyons \cite{AL} graphings of subexponential
growth are always hyperfinite.

\subsection{Graphings as graph limits}
Let $\bG=\gseq\subset\grd$ be weakly convergent graph sequence as in the
Introduction. Let $(\cG,X,\mu)$ be a graphing. If
$\alpha\in U^{r,d}$ then let
$T(\cG,\alpha)$ be the set of points $p\in X$ such that 
 the ball $B_r(p)\subset\cG_p$ is rooted isomorphic
to $\alpha$. Clearly, $T(\cG,\alpha)$ is a measurable set. We say that
$\bG$ converges to $\cG$ if for any $r\geq 1$ and $\alpha\in U^{r,d}$
$$p_{\bG}(\alpha):=\limn p_{G_n}(\alpha)=\mu(T(\cG,\alpha))\,.$$
In \cite{Elek1} we proved that any weakly convergent graph sequence
admits such {\bf limit graphings}. There is however an other
even more natural limit object
for weakly convergent graph sequences constructed by Benjamini and
Schramm \cite{BS}. Let $\gd$ be the set of all countable
connected rooted graphs (up to rooted isomorphism) with uniform
vertex degree bound $d$. For each $\alpha\in U^{r,d}$ we associate
a closed-open set $R(\alpha)$, the set of elements $G\in \gd$ such that
$B_r(x)\cong\alpha$, where $x$ is the root of $G$. Then $\gd$ is
a metrizable, compact space. Now let $\bG=\gseq$ be a weakly convergent
sequence in $\grd$. Then $$\hmu(R(\alpha)):=\limn p_{G_n}(\alpha)
=p_{\bG}(\alpha)$$
defines a measure $\hmu$ on $\gd$. This measure space
can be considered as the primary limit object
for weakly convergent graph sequences.

\noindent
Note that if $\bG=\gseq$ is a \Fo-sequence in the
Cayley-graph of a finitely generated
amenable group $\Gamma$ then the limit measure $\hmu$ is
concentrated on one single point in $\gd$ namely on the point
representing the Cayley-graph itself. 
In order to avoid this technical difficulty, in the following subsections
we introduce a combination of the limit graphing and the Benjamini-Schramm
construction.

\subsection{B-graphs}\label{Bgrap}
Let $B=\{0,1\}^\bN$ be the Bernoulli space of $0-1$-sequences with the
standard product measure $\nu$.
A rooted $B$-graph is a rooted connected graph $G$ equipped with a function
$\tau_G:V(G)\to B$. We say that the rooted $B$-graphs $G$ and $H$
are isomorphic if there exists a rooted graph isomorphism $\psi:G\to H$
such that $\tau_H(\psi(x))=\tau_G(x)$ for any $x\in V(G)$.
Let $\bgd$ be the set of such isomorphism classes of countable rooted
$B$-graphs with vertex degree bound $d$. 
Let $\alpha\in U^{r,d}$ and consider a rooted $r$-ball $T$ representing
the class $\alpha$. Consider the product space $B(T)=B^{V(T)}$
with the product measure $\nu^{|V(T)|}=\nu_T$.
Note that the finite group of rooted
automorphisms $\mbox{Aut}(T)$ acts continuously
on $B(T)$ preserving the measure $\nu_T$. Let us consider
the quotient space $Q(T)=B(T)/\mbox{Aut}(T)$ and the natural projection
$\pi_T:B(T)\to Q(T)$. For a Borel-set $W\subseteq Q(T)$ let us define
the measure $\lambda$ by $\lambda(W)=\nu_T(\pi_T^{-1}(W))\,.$
Obviously if $T$ and $S$ are rooted isomorphic balls then $Q(T)$ and $Q(S)$
are naturally isomorphic. Hence we shall
denote the quotient space by $Q(\alpha)$.
Let $\beta\in U^{r+1,d},\alpha\in U^{r,d}$ such that the $r$-ball around
the root in $\beta$ is isomorphic to $\alpha$. Then we have a natural
projection $\pi_{\beta,\alpha}:Q(\beta)\to Q(\alpha)$. Indeed if $f\in B(T)$
for some rooted ball $T$ representing $\beta$ and the restriction of $f$
on the $r$-ball around the root is $g$, then the class of $f$ is mapped to 
the class of $g$. 
\begin{lemma}\label{l11}
If $W\subseteq Q(\alpha)$ is a Borel-set then
$$\lambda(\pi^{-1}_{\beta,\alpha}(W))=\lambda(W)\,.$$
\end{lemma}
\proof
Let $\pi_{T,S}:B(T)\to B(S)$ be the natural projection, where $S$ is the
$r$-ball around the root. Then
$$\pi_S\circ \pi_{T,S}=\pi_{\beta,\alpha}\circ \pi_T\,.$$
Also, since $\nu_T(\pi^{-1}_{T,S}(A))=\nu_S(A)$ for
any Borel-set $A\subseteq B(S)$\,,
$$\lambda(\pi^{-1}_{\beta,\alpha}(W))=\nu_T(\pi^{-1}_T\circ
\pi^{-1}_{\beta,\alpha}(W))\,\,\mbox{and}\,\,
\lambda(W)=\nu_S(\pi^{-1}_S(W))=\nu_T(\pi^{-1}_{T,S}\circ \pi^{-1}_S(W))\,.$$
Now the lemma follows. \qed

\vskip 0.2in
\noindent
Hence we have the compact spaces
$Q^r_d:=\bigcup_{\alpha\in U^{r,d}} Q(\alpha)$ and the projections
$$Q^1_d\stackrel{\pi^1}{\leftarrow} Q^2_d\stackrel{\pi^2}{\leftarrow}\dots
\,$$
where $\pi^r$ is defined as $\pi_{\beta,\alpha}$ on $Q(\beta)$.
It is easy to see that the elements of $\lim_{\leftarrow} Q^r_d$
are in a one-to-one correspondence with the rooted isomorphism classes
of the countable rooted $B$-graphs with vertex degree bound $d$.
Hence from now on we regard $\bgd$ as a compact metrizable space. Note
that the forgetting functor provides us a continuous map
$\cF:\bgd\to\gd$. Note that the forgetting functor maps a $B$-graph to its
underlying graph in $\gd$.

\noindent
Now let $\bG=\gseq$ be a weakly convergent graph sequence and
$\hmu$ be the limit measure on $\gd$.
\begin{proposition}
Let $\alpha\in U^{r,d}$ and $W\subseteq Q(\alpha)$ be a Borel-set.
We define the measure $\tmu$ by
$$\tmu(W):=\lambda(W) p_{\bG}(\alpha)\,.$$
Then $\tmu$ is a Borel-measure on $\bgd$ and $\cF_*(\tmu)=\hmu\,.$
\end{proposition}
\proof
Clearly, we define a measure $\tmu^r$ on $Q^r_d$ by
$$\tmu^r(\cup_{\alpha\in U^{r,d}} W_\alpha):=\sum_{\alpha\in U^{r,d}}
p_{\bG}(\alpha)\lambda(W_\alpha)\,.$$
We only need to prove that
$$\pi^r_*(\tmu^{r+1})=\tmu^r\,.$$
Let $U_\alpha^{r+1,d}$ be the 
set of classes such that the rooted $r$-ball around
the root is just $\alpha$. Then
\begin{itemize}
\item $\bigcup_{\alpha\in U^{r,d}} U_\alpha^{r+1,d}=U^{r+1,d}\,.$
\item $(\pi^r)^{-1}(Q(\alpha))=\bigcup_{\beta\in U_\alpha^{r+1,d}} Q(\beta)\,.$
\item $p_{\bG}(\alpha)=\sum_{\beta\in U_\alpha^{r+1,d}} p_{\bG}(\beta)\,.$
\end{itemize}
If $W\subseteq Q(\alpha)$ then 
$$(\pi^r)^{-1}(W)=\bigcup_{\beta\in U_\alpha^{r+1,d}}
\pi^{-1}_{\beta,\alpha}(W)\,.$$
Thus by Lemma \ref{l11}, $\tmu^{r+1}(\pi^{-1}_{\beta,\alpha}(W))=
p_{\bG}(\beta)\lambda(W)\,.$
Therefore $$\tmu^{r+1}((\pi^r)^{-1}(W))=\tmu^r(W)\,.$$
Consequently, $\tmu$ is a well-defined Borel-measure on $\bgd$.
Since $\hmu(R(\alpha))=\tmu(Q(\alpha))$, $\cF_*(\tmu)=\hmu\,.\,\,\qed$
\subsection{The canonical colouring of a $B$-graph} \label{canon}
Let us consider the triples $(p,q,n)$, where
$1\leq p \leq d$, $1\leq q \leq d$, $n\geq 1$. Let $G(V,E,\tau_G)$ be
a countable $B$-graph such that $\tau_G(x)\neq\tau_G(y)$ if $x\neq y$. These
$B$-graphs are called {\bf separated}. Now colour the edge $e=(x,y)\in E$
by $(p,q,n)$ if
\begin{itemize}
\item $\tau_G(x)<\tau_G(y)$ (in the lexicographic ordering of $\{0,1\}^\bN$)
 and $l_1<l_2<\dots<l_{\mbox{deg}(x)}$ are
the values of $\tau_G$ at the neighbours of $x$ and $\tau_G(y)=l_p$.
\item $m_1<m_2<\dots<m_{\mbox{deg}(y)}$ are
the values of $\tau_G$ at the neighbours of $y$
and $\tau_G(x)=m_q$.
\item $\tau_G(x)=\{a_1,a_2,\dots\}\in B$, $\tau_G(y)=\{b_1,b_2,\dots\}\in B$,
$a_1=b_1, a_2=b_2,\dots, a_{n-1}=b_{n-1}, a_n\neq b_n\,.$
\end{itemize}
\begin{lemma}
If $(a,b)\in E$, $(a,c)\in E$ then the colours of $(a,b)$ and $(a,c)$
are different.
\end{lemma}
\proof
If the colour of $(a,b)$ and $(a,c)$ are the same, then either
$\tau_G(b)>\tau_G(a),\tau_G(c)>\tau_G(a)$ or
$\tau_G(b)<\tau_G(a),\tau_G(c)<\tau_G(a)$. Hence by the definition
of the colouring $\tau_G(b)=\tau_G(c)$ leading to a contradiction. \qed

\vskip0.2in
\noindent
Now consider ${\bf O}_d\subset \bgd$, the Borel-set of separated $B$-graphs.
Clearly $\tmu({\bf O}_d)=1$.
The colouring construction above defines a canonical Borel $\f2$-action
on ${\bf O}_d$ as follows. Suppose that $\uz\in {\bf O}_d$ represents the rooted
$B$-graph $G$ with root $a\in V(G)$. Consider the free generators
of order two $\{s_\delta\}_{\delta\in I}$, where 
$$I=\{1,2,\dots,d\}\times \{1,2,\dots,d\} \times \bN\,.$$
Let $\alpha\in I$, $\alpha=(p,q,n)$. Then
\begin{itemize}
\item If there exists an edge $(a,b)\in E(G)$ coloured by $(p,q,n)$ then define
  $s_\alpha(\uz)=\underline{w}$, where $\underline{w}$ represents the same
$B$-graph as $\uz$, but with root $b$.
\item If there exists no edge $(a,b)\in E(G)$ coloured by
$(p,q,n)$ then let $s_\alpha(\uz)=\uz$.
\end{itemize}

 Observe that we constructed
an $\f2$-action on ${\bf O}_d$ such a way that if
$\uz\in {\bf O}_d$ represents a graph $G$ then the orbit graph of $\uz$ is
isomorphic to $G$, We call this action the canonical $\f2$-action
on the {\bf canonical limit object} $({\bf O}_d,\tmu)$.
In Corollary \ref{pres} we shall prove that the measure $\tmu$ is
invariant under the canonical action.
\subsection{Random $B$-colourings of convergent graph sequences}
Let $\bG=\gseq\subset\grd$ be a weakly convergent sequence
of graphs. Let $\Omega=B^{\cup^\infty_{n=1} V(G_n)}$ be the
space of $B$-valued functions $\kappa$ on the vertices of
the graph sequence. We equip $\Omega$ with the standard product
measure $\nu_\Omega$.

\noindent
Now let $\alpha\in U^{r,d}$ and let $Q(\alpha)$ be the quotient
space as in Subsection \ref{Bgrap}. Let $U\subseteq Q(\alpha)$ be
a Borel-subset, $\kappa\in\Omega$ and $T(G_n,\kappa,U)$
be the set of vertices $p\in V(G_n)$ such that
\begin{itemize}
\item $p\in T(G_n,\alpha)\,.$
\item $\kappa_{\mid B_r(p)}\in U\,.$
\end{itemize}
\begin{proposition} \label{p19}
For any Borel-set $U\in Q(\alpha)$,
\begin{equation} \label{dodo}
\limn \frac {|T(G_n,\kappa,U)|}{|V(G_n)|}=\lambda(U)p_{\bG}(\alpha)=\tmu(U)
\end{equation}
holds for almost all $\kappa\in\Omega$.
\end{proposition}
\proof
We may suppose that $p_{\bG}(\alpha)\neq 0$, since if 
$p_{\bG}(\alpha)=0$ then both sides of the equation (\ref{dodo}) vanish.
Let $x\in T(G_n,\alpha)\,.$ Then we define $A_x^U\subset \Omega$ by
$$A_x^U:=\{\kappa\in\Omega\,\mid\, x\in T(G_n,\kappa,U)\}\,.$$
Clearly, $\nu_\Omega(A_x^U)=\lambda(U)$. Note however that if
$x\neq y\in T(G_n,\alpha)$ then $A_x^U$ and $A_y^U$ might not be
independent subsets. On the other hand, if
$x\in T(G_n,\alpha)$, $y\in T(G_m,\alpha)$ and $n\neq m$ then
$A_x^U$ and $A_y^U$ are independent. Also, if $S\subset T(G_n,\alpha)$,
$S=\{x_1,x_2,\dots,x_k\}$ and $d_{G_n}(x_i,x_j)>2r$ if $i\neq j$ then
$A_{x_1}^U, A_{x_2}^U,\dots, A_{x_k}^U$ are jointly independent.
\begin{lemma}
There exists a natural number $l>0$
(depending on $r$ and $d$) and a partition
$\cup^l_{i=1} B^n_i=T(G_n,\alpha)$ for any $n\geq 1$
such that if $x\neq y\in B^n_i$ then
$d_{G_n}(x,y)>2r\,.$
\end{lemma}
\proof
Let $H_n$ be a graph with vertex set $V(G_n)$.
Let $(x,y)\in E(H_n)$ if and only if $d_{G_n}(x,y)\leq 2r\,.$
Then $\mbox{deg}(x)\leq d^{r+1}$ for any $x\in V(H_n)$. Let
$l=d^{r+1}+1$ then $H_n$ is vertex-colorable by the colours
$c_1,c_2,\dots,c_l$.
Let $B_i^n$ be the set of vertices coloured by $c_i$.\qed

\vskip0.2in
\noindent
To conclude the proof of Proposition \ref{p19}, let us fix $q\geq 2$. Let
$B^n_{i_1}, B^n_{i_2}, \dots, B^n_{i_{n,q}}$ be those elements
of the partition of the previous lemma such that
$$\frac{|B^n_{i_j}|}{|V(G_n)|}>\frac{2^{-q}}{l}\,.$$
Then by the law of large numbers, for almost all $\kappa\in\Omega$
$$\limn\frac{|T(G_n,\kappa,U)\cap
  B^n_{i_j}|}{|B^n_{i_j}|}=\lambda(U)\,,$$
for any choice of $i_j$. An easy calculation shows that for the
same $\kappa$
\begin{equation} \label{dodo2}
\limn 
\frac{|T(G_n,\kappa,U)\cap
(\bigcup_{j=1}^{i_{n,q}}  B^n_{i_j})|}
{|\bigcup_{j=1}^{i_{n,q}}B^n_{i_j}|}=\lambda(U)\,.
\end{equation}
On the other hand,
$$\frac{|T(G_n,\alpha)\backslash \bigcup_{j=1}^{i_{n,q}}B^n_{i_j}|}
{|V(G_n)|}\leq 2^{-q}$$
and $\limn \frac{|T(G_n,\alpha)|}{|V(G_n)|}=p_{\bG}(\alpha)\,.$
Hence letting $q\to\infty$, (\ref{dodo}) follows.\qed
\subsection{Generic elements}
For any $\alpha\in U^{r,d}$ let us choose closed-open sets
$\{U^k_\alpha\}^\infty_{k=1}$ such that they form a Boolean-algebra
and generate all the Borel-sets in $Q(\alpha)$.
We call $\kappa\in\Omega$ {\bf generic} if for any $k\geq 1$ and
$\alpha\in U^{r,d}$
$$\limn \frac{|T(G_n\kappa, U^k_\alpha)|}{|V(G_n)|}=\lambda(U^k_\alpha)
p_{\bG}(\alpha)\,$$
and for any $n\geq 1$, $\kappa(p)\neq \kappa(q)$ if $p\neq q\in V(G_n)$.
By Proposition \ref{p19}, almost all $\kappa\in\Omega$ are generic.

\section{Graph sequences and ultraproducts}
\subsection{Basic notions}
In this section we briefly recall some of the basic notions on the
ultraproducts of finite sets \cite{ESZ}.
Let $\{X_i\}^\infty_{i=1}$ be finite sets, $|X_i|\to\infty$.
Let $\omega$ be a non-principal
ultrafilter and $\limo:l^\infty(\bN)\to\bR$ be the corresponding 
ultralimit. 
The ultraproduct of the sets $X_i$ is defined as follows.

\noindent
Let $\widetilde{X}=\prod^\infty_{i=1}X_i$. We say that
$\widetilde{p}=\{p_i\}^\infty_{i=1}, \widetilde{q}=\{q_i\}^\infty_{i=1}\in
\widetilde{X}$ are equivalent, 
$\widetilde{p}\sim\widetilde{q}$, if
$$\{i\in \bN\mid p_i=q_i\}\in \omega\,.$$
We shall denote the equivalence class of $\{p_i\}^\infty_{i=1}$ by
$[\{p_i\}^\infty_{i=1}]\,.$
Define $\bX:=\widetilde{X}/\sim$.
Now let $\cR(X_i)$ denote the Boolean algebra of subsets of $X_i$, with the
normalised measure $\mu_i(A)=\frac{|A|}{|X_i|}\,.$
Then let $\widetilde{\cR}=\prod^\infty_{i=1}\cR(X_i)$ and
$\cR=\widetilde{\cR}/I$, where $I$ is the ideal of elements
$\{A_i\}^\infty_{i=1}$
such that 
$\{i\in \bN\mid A_i=\emptyset\}\in \omega\,.$
It is important to note
 that the elements of $\cR$ can be identified with certain subsets
of $\bX$: If 
$$\bp=[\{p_i\}^\infty_{i=1}]\in \bX\,\,\mbox{and}\,\, \bA=
[\{A_i\}^\infty_{i=1}]\in \cR$$
then $\bp\in \bA$ if 
$\{i\in \bN\mid p_i\in A_i\}\in \omega\,.$
One can easily see that $\cR$ is a Boolean-algebra on $\bX$.
Now let $\bmu(\bA)=\limo \mu_i(A_i)\,.$
Then $\bmu:\cR\to\bR$ is
a finitely additive probability measure.
We call $\bbN\subseteq \bX$ a nullset if for any $\e>0$ there exists
$\bA_\e\in\cR$ such that
$\bbN\subset \bA_\e$ and $\mu(\bA_\e)\leq\e$.
We call $\bB\subset \bX$ {\bf measurable} if there exists
$\hat{\bB}\subset\cR$ such that $\bB\triangle \hat{\bB}$ is a nullset.
The measurable sets form a $\sigma$-algebra $\cB$ and
$\bmu(\bB)=\bmu(\hat{\bB})$
defines a probability measure on $\cB$.
\subsection{The ultraproduct of $B$-valued functions}
Let $\bG=\gseq\subset \grd$ be a weakly convergent sequence of graphs.
We shall denote by $\xg$ the ultraproduct of the vertex sets
$\{V(G_n)\}^\infty_{n=1}$.
Now consider an element $\kappa\in\Omega=B^{\cup^\infty_{n=1}V(G_n)}\,.$
We define the $B$-valued function $F_\kappa$ on $\xg$ the following way.
Let $\bp=[\{p_n\}^\infty_{n=1}]$ then $F_\kappa(\bp):=\limo \kappa(p_i)$.
Note that if $\{b_n\}^\infty_{n=1}\subset B$ is a sequence of elements
of the Bernoulli product space then $\limo b_n=b$ is the unique element
of $B$ such that for any neighbourhood $b\in U\subseteq B$
$$\{n\in\bN\,\mid\, b_n\in U\}\in \omega\,.$$
\begin{lemma}
$F_\kappa$ is a measurable $B$-valued function on $\xg$.
\end{lemma}
\proof
Let $O_{x_1,x_2,\dots,x_n}$ be the
basic closed-open set in $B$, where $x_i\in\{0,1\}$ and
$b\in O_{x_1,x_2,\dots,x_n}$ if $b(i)=x_i$. 
It is enough to prove that $F_{\kappa}^{-1}(O_{x_1,x_2,\dots,x_n})\in\cR\,.$
Let
$$ O^i_{x_1,x_2,\dots,x_n}:=\{p_n\in V(G_n)\,\mid\,\kappa(p_i)\in
O_{x_1,x_2,\dots,x_n}\}\,.$$
Since $ O^i_{x_1,x_2,\dots,x_n}$ is an closed-open set
$[\{O^i_{x_1,x_2,\dots,x_n}\}^\infty_{i=1}]=
F_{\kappa}^{-1}(O_{x_1,x_2,\dots,x_n})$. Thus our lemma follows. \qed
\subsection{The canonical action on the ultraproduct space}
Let $\bG=\gseq\subset\grd$ and $\xg$ be as in the previous subsections and
let $\kappa\in\Omega$ be a fixed generic element.
Then $\kappa$ determines a separated $B$-function on each vertex space
$V(G_n)$. Now let $\beta\in\{1,2\dots,d\}\times\{1,2,\dots,d\}\times\bN$
and let $S^i_\beta:V(G_i)\to V(G_i)$ be the
bijection $S^i_\beta(p)=s_\beta(p)\,$(as in Subsection \ref{canon}).
The ultraproduct of $\{S^i_\beta\}^\infty_{i=1}$ is defined
the following way:
$$\bS_\beta([\{p_i\}^\infty_{i=1}])=[\{S^i_\beta(p_i)\}^\infty_{i=1}]\,.$$
Then $\bS_\beta$ is a measure-preserving bijection on the ultraproduct space
$\xg$. Indeed if $\bA=[\{A_i\}^\infty_{i=1}]\in\cR$ then
$\bS_\beta(\bA)=[\{S^i_\beta(A_i)\}^\infty_{i=1}]\,.$ 
Clearly $\bS_\beta^2=\mbox{Id}$, hence
we defined a measure-preserving action of $\f2$ on $\xg$.
\begin{lemma}
Each component $\cG_{\bp}$ of the graphing $\cG$ induced by the action
above has vertex degree bound $d$.
\end{lemma}
\proof Let $\bp=[\{p_i\}^\infty_{i=1}]\in\xg$. If
$\bS_\beta(\bp)\neq \bp$ then $s_\beta(p_i)\neq p_i$ for
$\omega$-almost all $i\in \bN$. Therefore if
$\bS_{\beta_1}, \bS_{\beta_2},\dots,\bS_{\beta_{d+1}}$ are bijections
such that $\bS_{\beta_j}(\bp)\neq\bp$ then $s_{\beta_j}(p_i)\neq p_i$
for $\omega$-almost all $i\in\bN$ and $1\leq j \leq d+1$. This leads
to a contradiction. \qed

\vskip0.2in
\noindent
By the previous lemma each graph $\cG_{\bp}$ is a rooted $B$-graph of vertex
degree
bound $d$, where $\bp$ is the root and the $B$-colouring on the vertices of
 $\cG_{\bp}$ is induced by $F_\kappa$.
Consequently, we have a canonical map $\rho:\xg\to\bgd$ (depending on the fixed
generic element $\kappa$ of course) such that for each $\bp$,
$\rho(\bp)$ is the rooted $B$-graph representing the component $\cG_{\bp}$.
\subsection{The canonical map preserves the measure} \label{preserve}
The goal of this subsection is to prove the main technical tool of our
paper.
\begin{proposition} \label{tool}
$\rho:(\xg,\mu_{\bG})\to (\bgd,\tmu)$ is a measure-preserving map.
\end{proposition}
\proof We need to prove that for any Borel-set $W\subseteq\bgd$,
$\rho^{-1}(W)$ is a measurable set in $\xg$ and
$\mu_{\bG}(\rho^{-1}(W))=\tmu(W)\,.$
\begin{lemma} \label{l33}
Suppose that the $r$-neighborhood of $\bp=[\{p_i\}^\infty_{i=1}]\in\xg$
represents $\alpha\in U^{r,d}\,.$ Then for $\omega$-almost all $i\in\bN$ the
$r$-neighbourhood of $p_i\in V(G_i)$ represents $\alpha$ as well.
\end{lemma}
\proof
Let $\bq\in B_r(\bp)$. Then there exists a path $\bp^0,\bp^1,\dots,\bp^r$
in the graph $\cG_{\bp}$ such that $\bp^0=\bp$, $\bp^r=\bq$. Therefore
for $\omega$-almost all $i\in \bN$, $(p^k_i,p^{k+1}_i)\in E(G_i)$
thus if $\bq=[\{q_i\}^\infty_{i=1}]$ then $q_i\in B_r(p_i)$ for
$\omega$-almost $i\in\bN$. Obviously if $\bq$ and $\bq'$ are vertices
in $B_r(\bp)$ then $(\bq,\bq')\in E(\cG_{\bp})$ if and only if
$(q_i,q_i')\in E(G_i)$ for $\omega$-almost all $ i \in\bN$. Also,
if $\mbox{deg}(\bq)=k$ then $\mbox{deg}(q_i)=k$ for
 $\omega$-almost all $ i \in\bN$.
This shows that $B_r(p_i)\cong B_r(\bp)$ for $\omega$-almost all $i\in\bN$.
\qed
\begin{lemma}
\label{l34}
Let $U\subseteq Q(\alpha)$ be a closed-open subset. Then
$F_{\kappa\mid B_r(\bp)}\in U$ if and only if
$\kappa_{\mid B_r(p_i)}\in U$ for $\omega$-almost all $i\in\bN$,
where $\kappa_{\mid B_r(p_i)}$ denotes the restriction of $\kappa$
onto the set $ B_r(p_i)$.
\end{lemma}
Observe that ${F_\kappa}_{\mid B_r(\bp)}=\limo \kappa_{\mid B_r(p_i)}$. Note
that the ultralimit of a sequence in a compact metric space is in the
closure of the sequence. Therefore the lemma easily follows. \qed.

\vskip0.2in
\noindent
By Lemma \ref{l34}, if $U\subseteq Q(\alpha)$ be a closed-open subset
$$\bmu(\{\bp\in\xg\,\mid\,B_r(\bp)\cong\alpha\,\,\mbox{and}\,\,
F^\kappa_{\mid\,B_r(\bp)}\in U\})=\limo\frac{|T(G_i,\kappa, U)|}{|V(G_i)|}\,.$$
Since $\kappa$ is a generic element in $\Omega$,
$\mu_{\bG}(\rho^{-1}(U^k_\alpha))=\tmu(U^k_\alpha)$ for
any $\alpha\in U^{r,d}$ and $k\geq 1$.
Since $\{U^k_\alpha\}^\infty_{k=1}$ is a generating Boolean-algebra
$\mu_{\bG}(\rho^{-1}(W))=\tmu(W)$ holds
for any Borel-set $W\subseteq \bgd$.\qed
\begin{corollary}{} \label{pres}
\vskip0.1in
\noindent
\begin{description}
\item{(a)} For almost all $\bp\in\xg$, $\cG_{\bp}$ is a separated
$B$-graph.
\item{(b)} The $\f2$-action on ${\bf O}_d\subseteq \bgd$ preserves the measure.
\end{description}
\end{corollary}
\proof (a) follows from the fact that $\mu_{\bG}(\rho^{-1}({\bf O}_d))=1$. 
On the other
hand $\rho$ commutes with the $\f2$-action, that implies (b).\qed
\subsection{The ultraproduct of finite graphs}
The goal of this subsection is to prove some auxiliary lemmas that shall be 
used in the proof of our main theorem. Let $\bG=\gseq$ be a weakly
convergent sequence of finite graphs with vertex degree bound $d$. Let $\xg$
be the ultraproduct of their vertex sets and $\cG_{\bX}$ be the graphing
constructed in the previous subsection.
Then we have two notions for measure space of edge-sets. The first is the
one constructed in Subsection \ref{hyperfinite}. On the other hand,
similarly to the ultraproduct of the vertex sets we can also define
the ultraproducts of edge sets with normalised measure $\mu_{\bE}$, 
$$\mu_{\bE}(L)=\limo \frac {|E(L_n)|}{|V(G_n)|}\,,$$
where $L=[\{L_n\}^\infty_{n=1}]$, $L_n\subseteq E(G_n)$\,. Again the
ultraproduct sets $L=[\{L_n\}^\infty_{n=1}]$ form a Boolean-algebra
$\cR_{\bE}$ and we can define the $\sigma$-algebra of measurable edge-sets
by $\cM_{\bE}$ as well. 

\noindent
It is easy to see that the two measure spaces above coincide.
If $\bA\in\cM_{\bE}$, then let $V(\bA)$ be the set of 
points $\bp$ in $\xg$ for which there exists $\bq\in\xg$, $(\bp,\bq)\in \bA$.
Clearly if $\bA\in\cR_{\bE}$ then $V(\bA)\in\cR$ and if $\bbbN$ is a 
nullset of edges
then $V(\bbbN)$ is a nullset of vertices. Consequently if $A\in\cM_{\bE}$ then
$V(\bA)\in\cM$. Note that we can regard the elements of $\cM_{\bE}
$ as measurable
subgraphs of $\cG_{\bX}$.
\begin{lemma}
\label{appro1}
Let $\bH_F\in \cM_{\bE}$ be a subgraph
such that all of its components are isomorphic to
a finite simple graph $F$. Then for any $\gamma>0$ there exists 
$\bS_F\subset \bH_F$ such that
\begin{itemize}
\item
$\bS_F\in\cR_{\bE}$
\item All the components of $\bS_F$ are isomorphic to $F$.
\item $\mu_{\bG}(V(\bH_F\backslash \bS_F))<\gamma\,.$
\end{itemize}
\end{lemma}
\proof
Let $\bH_F'\in\cR_{\bE}$ be a subgraph such that 
$\mu_{\bE}(\bH_F\triangle \bH'_F)=0\,.$
Then $V(\bH_F\triangle \bH'_F)$ is a nullset in $\xg$.
Consequently, $\bQ=\mbox{Orb}(V(\bH_F\triangle \bH'_F))$ is
still a nullset in $\xg$, where
$$\bQ=\cup_{\bp\in V(\bH_F\triangle \bH'_F)} V(\cG_\bp)$$
is the union of the orbits of the vertices in $\bH_F\triangle \bH'_F$.
Let $\bQ\subset\bK\in\cR$, $\mu_{\bG}(\bK)\leq \gamma$.
Consider the subset $V(\bH')\backslash \bK\in\cR$. Then
the $r$-neighbourhood of $V(\bH')\backslash \bK$, 
$B_r(V(\bH')\backslash \bK)$ is
also an element of $\cR$ for any $r\geq 1$.
Indeed, if $V(\bH')\backslash \bK=[\{A_i\}^\infty_{i=1}]$
then
$$B_r(V(\bH')\backslash \bK)= [\{B_r(A_i)\}^\infty_{i=1}]\,.$$
Let $r>\mbox{diam}\,(F)$ and $\bS_F$ be the spanned subgraph
of $B_r(V(\bH')\backslash \bK)$ in $\bH'$. Then clearly $\bS_F\in\cR_{\bE}$.
Also, $\bS_F$ does not contain any vertex of
$V(\bH'\triangle \bH)\,.$ Thus $V(\bS_F)\subseteq V(\bH)$ and
if $\bp\in V(\bS_F)$ then the component
of $\bp$ in $\bS_F$ is just the component of $\bp$ in $\bH$. Clearly,
$\mu_{\bG}(V(\bH\backslash \bS_F))<\gamma$
thus our lemma follows.\,\,\qed
\begin{lemma}
\label{appro2}
Let $\bS=[\{S_n\}^\infty_{n=1}]\in \cR_{\bE}$ be a subgraph such that
all of its components are isomorphic to the finite simple graph $F$.
Then for $\omega$-almost all $n$ each component of $S_n$ is isomorphic to
$F$\,.
\end{lemma}
\proof
We prove the lemma by contradiction.
Suppose that there exists $T\in\omega$ such that for any $n\in T$ there exists
$p_n\in V(S_n)$ such that the component
of $S_n$ containing $ p_n$ is not isomorphic to $F$.

\vskip0.2in
\noindent
\underline{Case $1$}: If for $\omega$-almost all elements of $T$ 
there exists $q_n\in
V(S_n)$ such that $d(p_n,q_n)=r>\mbox{diam} (F)$ then there exists
$(\bp,\bq)\in\bS$ such that
$d_\bS(\bp,\bq)>\mbox{diam} (F)$, leading to a contradiction.
\vskip0.2in
\noindent
\underline{Case $2$}: 
If for $\omega$-almost all elements of $T$ the component containing
$p_n$ has diameter less than $3r$, then for $\omega$-almost all elements of
$T$
the component containing $p_n$ is isomorphic to the same finite graph $G$,
where $G$ is not isomorphic to $F$. Then there exists $\bp\in V(\bS)$
such that the component of $\bS$ containing $p$ is isomorphic to $G$.
This also leads to a contradiction. \qed

\section{The proof of Theorem \ref{main}}
Let $\e>0$ and $\bG=\gseq\subset\grf$ be a weakly convergent sequence
of graphs. Also let $(\bgd,\tmu)$ be the canonical limit object
as in Subsection \ref{canon}.
Let $\bgf$ be the $\f2$-invariant subspace of graphs $G$ in $\bgd$ satisfying
$|B_r(x)|\leq f(r)$ for all $x\in V(G)$.
Let ${\bf O}_f=\bgf\cap {\bf O}_d$. Then $\bO_f$ is also $\f2$-invariant and
$\tmu(\bgd\backslash \bO_f)=0\,.$
Consider the induced graphing $(\cG,\bO_f,\tmu)$. Since all the component
graphs are of subexponential growth, by the theorem of Adams and Lyons
\cite{AL}, this graphing is hyperfinite.
Therefore there exists a $K>0$ and a component-finite subgraphing
$\cH\subset\cG$ such that
\begin{itemize}
\item $\mu_E(E(\cG)\backslash E(\cH))\leq \frac{\e}{2}\,.$
\item $\bO_f=\bigcup_{H,\,|V(H)|\leq K} \cH_H$, where
$\cH_H$ is the set of points in $\bO_f$ contained in a component of
$\cH$ isomorphic to $H$.
\item $E(\cH)=\bigcup_{H,\, |V(H)|\leq K} E(\cH_H)$
\end{itemize}
Let $c_H=\tmu(\cH_H)$. Now suppose that our theorem does not hold.
Therefore there exists a subsequence $\{G_{n_i}\}^\infty_{i=1}$ such that
one can not remove $\e E(G_{n_i})$ edges from any $G_{n_i}$ to satisfy
condition (1) of our Theorem with the extra condition that
$$\left| \frac{|V^{n_i}_H|}{|V(G_{n_i})|}-c_H\right| <\delta,$$
for any finite simple graph $H$, $|V(H)|\leq K$.
Let $\xg$ be the ultraproduct of the graphs $\{G_{n_i}\}^\infty_{n=1}$.
Note that the canonical limit objects of the sequence
$\{G_{n_i}\}^\infty_{i=1}$
and of $\gseq$ are the same.
Therefore we have a measure-preserving map
$$\rho:(\bO_{\bG},\mu_{\bG})\to (\bO_f,\tmu)\,,$$
where $\bO_{\bG}$ is the set of elements $\bp\in\xg$ such that
$\cG_{\bp}$ is separated. Note that
\begin{itemize}
\item $\mu_{\bG}(\xg\backslash \bO_{\bG})=0\,.$
\item $\rho$ commutes with the canonical $\f2$-actions.
\item$\rho$ preserves the isomorphism type of the orbit graphs.
\end{itemize}
Clearly, $\rho$ extends to a measure-preserving map
$$\hat{\rho}:(E(\bO_{\bG}),\mu_{\bE})\to (E(\bO_f),\tilde{\mu}_E)\,,$$
where $\mu_{\bE}$ and $\tilde{\mu}_E$ denote the induced measures on
the edge-sets.
Now fix a constant $\gamma>0$. Let $\bA_H=\hat{\rho}^{-1}(\cH_H)$.
Then $\{\bA_H\}_{H,|V(H)|\leq K}$ are component-finite subgraphings
and all the components of $\bA_H$ are isomorphic to $H$.
Observe that $\mu_{\bG}(V(\bA_H))=c_H\,.$ Now first apply Lemma \ref{appro1}
to obtain subgraphings $\bS_H\subset\bA_H$ such that
$\mu_{\bG}(V(\bA_H\backslash \bS_H))<\gamma$ for each $H$. Then we apply
Lemma \ref{appro2} to obtain the graphs $\{S^H_{n_i}\}^\infty_{i=1}$ 
for each $H$ such that
\begin{itemize}
\item
$S^H_{n_i}\subset G_{n_i}$
\item
All the components of $S^H_{n_i}$ are isomorphic to $H$.
\item
$\limo \left| \frac{|V(S^H_{n_i})|}{|V(G_{n_i})|}-c_H\right|<\gamma\,.$
\end{itemize}
Thus for $\omega$-almost all $i\in \bN$
\begin{itemize}
\item
$\left|\frac{|V(S^H_{n_i})|}{|V(G_{n_i})|}-c_H\right|< 2\gamma$
\item
$\left|\frac { E(G_{n_i})\backslash \bigcup_{ H, |V(H)|\leq K} E(S^H_{n_i})}
{|V(G_{n_i})|}\right|<2d\gamma g_K +\frac{\e}{2}\,$
where $g_K$ is the number of graphs having vertices not greater than $K$.
\end{itemize}
Since $\gamma$ can be chosen arbitrarily we are in contradiction with
 our assumption on the graphs $\{G_{n_i}\}^\infty_{i=1}$.\,\,\qed

\vskip 0.2in
\noindent
{\bf Remark :} In \cite{Sch}, Schramm proved that a graph sequence
$\{G_n\}^\infty_{n=1}$ is hyperfinite if and only if its unimodular limit 
measure is hyperfinite. This idea
was used in \cite{BSS} to show that planarity is a testable property for
bounded degree graphs. If we prove that the canonical limit of a
hyperfinite graph sequence is always hyperfinite, then we can extend the
results of our paper to arbitrary hyperfinite classes. This is subject of
ongoing research \cite{Elek3}. In \cite{CSS}, the authors studied
hereditary hyperfinite classes (see Corollary 3.2 of their paper).
A graph class is hereditary if it is closed under vertex removal.
Thus planar graphs of bounded degree $d$ and $\grf$ are both hereditary
hyperfinite classes. The main result of \cite{CSS} is that hereditary
properties are testable in hyperfinite classes. It means that a tester
accepts the graph if it has the property and rejects the graph with probability
at least $(1-\e)$ if the graph is $\e$-far from the property in
edit-distance. It would be interesting to see whether the edit-distance from
a hereditary property is testable in a hereditary hyperfinite graph class.

\section{Testing union-closed monotone graph properties}
\subsection{Edit-distance from a union-closed monotone graph property}
\begin{theorem}\label{theo2}
Let $\cP$ be a union-closed monotone graph property as in the Introduction.
Then $\zeta(G)=d_e(G,\cP)$ is a continuous graph parameter
on $\grf$.
\end{theorem}
\proof We define the {\bf normal distance} from a union-closed monotone class 
by
$$d_n(G,\cP)=\inf_{H\subset G, H\in\cP} d_e(G,H)\,.$$
\begin{lemma}\label{normal=edit}
$d_n(G,\cP)=d_e(G,\cP)$
\end{lemma}
\proof Clearly, $d_e(G,\cP)\leq d_n(G,\cP)$. Now let
$J\in\cP$, $V(J)=V(G)$. 
Then the spanning graph $J\cap G$ has also property $\cP$
and $d_e(G,J\cap G)\leq d_e(G,J)$. Therefore $d_e(G,\cP)\geq d_n(G,\cP)$.
\,\,\,\qed
\vskip0.2in
\noindent
Now we prove a simple continuity lemma.
\begin{lemma} \label{conti}
If $G'\subseteq G$, $d_e(G,G')\leq\delta$ then
$|d_n(G',\cP)-d_n(G,\cP)|\leq \delta$.
\end{lemma}
\proof Let $H'\subseteq G'$, $H'\in \cP$. Then
$d_e(G,H')\leq d_e(G',H')+\delta$. Consequently,
$d_n(G,\cP)\leq d_n(G',\cP)+\delta$.
Now let $H\subseteq G$, $H\in\cP$. Then $H\cap G'\in\cP$. Since
$d_e(G',H\cap G')\leq d_e(G,H)$, we obtain that $ d_n(G',\cP)\leq
d_n(G,\cP)\,.$\qed
\begin{lemma} \label{copy}
Let $A_1, A_2, A_3,\dots, A_l$ be finite simple graphs. Suppose that
the graph $A$ consists of $m_1$ disjoint copies of $A_1$ and $m_2$
disjoint copies of $A_2$ \dots and $m_l$ disjoint copies of $A_l$.
That is $|V(A)|=\sum_{i=1}^l m_i|V(A_i)|\,.$
Then
$$d_n(A,\cP)=\sum^l_{i=1} w_i d_n(A_i,\cP)\,,$$
where $w_i=\frac{ m_i|V(A_i)|}{\sum^l_{i=1} m_i |V(A_i)|}\,.$
\end{lemma}
\proof
Let $B\subset A$ be the closest subgraph in $\cP$. Then $B\cap A^j_i$ is the
closest subgraph in $\cP$ in each copy of $A_i$.
Hence
$$d_n(A,\cP)=\frac{|E(A\backslash B)|}{\sum^l_{i=1} m_i |V(A_i)|}=
\frac{\sum^l_{i=1} m_i E(A_i\backslash B)}{\sum^l_{i=1} m_i |V(A_i)|}=$$
$$
=\frac{\sum^l_{i=1} m_i d_n(A_i,\cP)|V(A_i)|}{\sum^l_{i=1} m_i |V(A_i)|}\,
\quad\qed$$
Now let $\bG=\gseq\in\grf$ be a weakly convergent
graph sequence and $\e>0$. Consider the graphs $G'_n$ in Theorem \ref{main}.
Then by Lemma \ref{conti}, $|d_n(G_n,\cP)-d(G_n',\cP)|<\e$.
Let $s^n_H$ be the number of components in $G'_n$ isomorphic to $H$.
By Lemma \ref{copy},
$$d_n(G'_n,\cP)=\sum_{H, |H|\leq K}
\frac{s^n_H |V(H)|d_n(H,\cP)}{|V(G_n)|}\,.$$
By Theorem \ref{main}, 
$$\limn \frac{s^n_H|V(H)|}{|V(G_n)|}=c_H\,.$$
Therefore $\limn d_n(G_n',\cP)= \sum_{H, |H|\leq K} c_H d_n(H,\cP)\,.$
Hence if $n,m$ are large enough then $|d_n(G_n,\cP)-d_n(G_m,\cP)|<3\e$.
Consequently, $\limn d_n(G_n,\cP)$ exists.\qed
\subsection{Testability versus continuity}

Let $\zeta:\grf\to\bR$ be a continuous graph parameter.
Let $\e>0$ be a real constant
and $N(\zeta,\e)>0$, $r(\zeta,\e)>0$, $k(\zeta,\e)>0$
be integer numbers. An $(\e,N,r,k)$-{\bf random sampling} is the
following process. For a graph $G\in \grf$, $|V(G)|\geq N(\zeta,\e)$ we
randomly pick $k(\zeta,\e)$ vertices of $G$. Then by examining
the $r(\zeta,\e)$-neighbourhood of the chosen vertices we obtain
an empirical distribution
$$Y:\bigcup_{s\leq r} U^{s,d}\to\bR\,.$$
A $(\zeta,\e)$-{\bf tester} is an algorithm $T$ which takes
the empirical distribution $Y$ as an input and calculates the
real number $T(Y)$. We say that $\zeta$ is {\bf testable}
if for any $\e>0$ there exist constants 
$N(\zeta,\e)>0$, $r(\zeta,\e)>0$, $k(\zeta,\e)>0$ and a $(\zeta,\e)$-tester
such that
$$\mbox{Prob}(|T(Y)-\zeta(G)|>\e)<\e\,.$$
In other words, the tester estimates the value of $\zeta$ on $G$ using
a random sampling and guarantees that the error shall be less than $\e$
with probability $1-\e$.
\begin{theorem} \label{test}
Any continuous graph parameter $\zeta$ on $\grf$ is testable.
\end{theorem}
\proof Since $\zeta$ is continuous on the compactification
of $\grf$, for any $\e>0$ there exist constants $r(\zeta,\e)>0$
and $\delta(\e,\zeta)>0$ such that
\begin{equation} \label{kulcs}
\mbox{If}\,\,|p_G(\alpha)-p_{G'}(\alpha)|<\delta\,\,\mbox{for all}\,\,
\alpha\in U^{s,d}, s\leq r
\,\,\mbox {then}\,\, |\zeta(G)-\zeta(G')|<\e\,.\end{equation}
Also, by the total boundedness of compact metric spaces, there exists
a finite family of graphs (depending on $\zeta$ and $\e$)
$\{G_1, G_2,\dots, G_t\}\subset\grf$ such that
 for any $G\in\grf$ there exists at least one $G_k$, $1\leq k \leq t$
such that 
$$|p_{G_k}(\alpha)-p_{G}(\alpha)|<\frac{\delta} {2}$$
for each $\alpha\in U^{s,d}, s\leq r$.
By the law of large numbers there exist constants $N(\zeta,\e)>0$ and
$k(\zeta,\e)>0$ such that if $Y$ is the empirical
distribution of an $(\e,N,r,k)$-random sampling then the probability
that there exists an $\alpha\in U^{s,d}$ for some $s\leq r$
satisfying
$$|p_G(\alpha)-Y(\alpha)|>\frac{\delta}{2}$$
is less than $\e$. Note that the sampling is taking place
 on the vertices of $G$
and $|V(G)|>N(\zeta,\e)$.

\noindent
The tester works as follows. First the sampler measures $Y$. Then the
algorithm compares the vector $\{Y(\alpha)\}_{\alpha\in U^{s,d},s\leq r}$
to a finite database containing the vectors
$$\{p_{G_1}(\alpha)\}_{\alpha\in U^{s,d},s\leq r},
\{p_{G_2}(\alpha)\}_{\alpha\in U^{s,d},s\leq r},\dots,
\{p_{G_t}(\alpha)\}_{\alpha\in U^{s,d},s\leq r}\,.$$
Now with probability at least $(1-\e)$ the algorithm finds $1\leq k\leq t$
such that
$|p_G(\alpha)-p_{G_k}(\alpha)|\leq\delta$ for any 
$\alpha\in U^{s,d}, s\leq r$. Then the output $T(Y)$ shall be $\zeta(G_k)$.
By (\ref{kulcs}) the probability that
$|\zeta(G)-T(Y)|>\e$ is less than $\e$.\qed

\section{Continuous parameters in $\grf$}\label{Ind}
\subsection{Integrated density of states}
Integrated density of states is a fundamental concept in mathematical
physics. Let us explain, how this notion is related to graph parameters.
Recall that the Laplacian
on the finite graphs $G$, $\Delta_G:l^2(V(G))\to l^2(V(G))$ is a 
positive, self-adjoint
operator defined by
$$\Delta_G(f)(x):=\mbox{deg}(x)f(x)-\sum_{(x,y)\in E(G)} f(y)\,.$$
For a finite dimensional self-adjoint linear operator $A:\bR^n\to \bR^n$
, the normalised spectral distribution of $A$ is given by 
$$N_{A}(\lambda):=\frac{s_{A}(\lambda)}{n}\,,$$
where $s_A(\lambda)$ is the number of eigenvalues of $A$ not greater than
$\lambda$ counted with multiplicities.
Therefore $N_{\Delta_G}(\lambda)$ is a graph parameter for every $\lambda\geq
0$. Now consider the $3$-dimensional lattice graph $\bZ^3$. The finite cubes
$C_n$ are the graphs induced on the sets $\{-n,-n+1,\dots,n-1,n\}^3$. It is
known for decades that for any  $\lambda\geq 0$ $\lim_{n\to\infty} 
N_{\Delta_{C_n}}(\lambda)$ exists and in fact the convergence is uniform in
$\lambda$. In other words, the integrated density of states exists in the
uniform sense.
The discovery of quasicrystals led to the study of certain infinite graphs 
that are not periodic as the lattice
graph. What sort of graphs are we talking about ? 

\noindent
Let $G$ be an infinite connected graph such that
$|B_r(x)|\leq f(r)$, for any $x\in V(G)$. We say that a
sequence of finite induced subgraphs $\{F_n\}^\infty_{n=1}$, 
$\cup^\infty_{n=1}F_n=G$ form
a \Fo-sequence if $\limn \frac{|\partial F_n|}{|V(F_n)|}=0\,,$
where
$$\partial F_n:=\{p\in V(G)\,\mid\, p\in V(F_n)\,\mbox{and there exists}
\, q\notin V(F_n) \,\,\mbox{such that}\,\,(p,q)\in E(G)\,\}$$
Note that subexponential growth implies that for any $x\in V(G)$,
$\{B_n(x)\}^\infty_{n=1}$ contains a \Fo-subsequence.

\noindent
We say that an infinite graph $G$ of subexponential growth
has {\bf uniform patch frequency} if all of its
\Fo-subgraph sequences are weakly convergent. Obviously, the lattices $\bZ^n$
are of uniform patch frequency, but there are plenty of aperiodic UPF graphs 
as well,
among them the graph of a Penrose tiling, or other Delone-systems \cite{LS}.

\noindent
Using ergodic theory, Lenz and Stollmann proved the existence of the 
integrated density of states in
the uniform sense for such Delone-systems \cite{LS}
and later we extended their results 
for all UPF graphs of subexponential growth \cite{Elek2}.
This last result can be interpreted the following way : 
$N_{\Delta_G}(\lambda)$ are continuous graph parameters in $\grf$ for
any $\lambda\geq 1$.

\vskip0.2in
In this subsection we apply our Theorem \ref{main} to extend the
theorem in \cite{Elek2}
for discrete Schr\"odinger operators with random potentials
 (as a general reference see the
lecture notes of Kirsch \cite{Kir}).
 Let $X$ be a random
variable taking finitely many real values
$\{r_1,r_2,\dots,r_m\}$. Let 
$$\mbox{Prob}(X=r_i)=p_i\,.$$
For the vertices $p$ of $G$ we consider independent random variables
$X_p$ with the same distribution as $X$. Let $\Omega^X_G$ be the
space of $\{r_1,r_2,\dots,r_m\}$-valued functions on $V(G)$ with the
product measure $\nu_G$. That is
$$\mu_G(\{\omega\,\mid\omega(x_1)=r_{i_1},
\omega(x_2)=r_{i_2}, \dots, \omega(x_k)=r_{i_k}\})=\prod^k_{j=1} p_{i_j}\,$$
for any $k$-tuple $(x_1,x_2,\dots,x_k)\subset V(G)$.
Thus for each $\omega\in\Omega^X_G$ we have a self-adjoint operator
$\Delta^\omega_G:l^2(V(G))\to l^2(V(G))$, given by
$$\Delta^\omega_G(f)(x)=\Delta_G(f)(x)+\omega(x)f(x)\,.$$
This operator is a discrete Schr\"odinger operator with random potential.
The following theorem is the extension of the main theorem of \cite{Elek2}
for such operators. Note that in
the case of Euclidean lattices a similar result was proved by Delyon and
Souillard \cite{DS}.
\begin{theorem} \label{int}
Let $f$ be a function of subexponential growth.
Let $G$ be an infinite connected graph such that
$|B_r(x)|\leq f(r)$, for any $x\in V(G)$
with UPF and let $\{G_n\}^\infty_{n=1}$ be a \Fo-sequence. Then for almost 
all $\omega\in\Omega^X_G$
$\{N_{\Delta^\omega_{G_n}}\}^\infty_{n=1}$ uniformly converges to
an integrated density of states $N^X_G$ that does not depend on $\omega$
That is the integrated density of state for such discrete Schr\"oedinger
operator with random potential is
non-random. (see also \cite{LV} and the references therein)
\end{theorem}
\proof Let $\e>0$ and $G_n'\subset G_n$ be the spanning graphs
as in Theorem \ref{main}. First we prove the analog of Lemma \ref{conti}.
\begin{lemma}\label{analoglemma}
For any $\omega\in\Omega^X_G$,
$$|N_{\Delta^\omega_{G_n}}(\lambda)-N_{\Delta^\omega_{G'_n}}(\lambda)
|\leq \e d\,,$$
for any $-\infty<\lambda <\infty\,,$ where
$d$ is the uniform bound on the degrees of $\{G_n\}^\infty_{n=1}$.
\end{lemma}
\proof
Observe that
\begin{equation}\label{rank}
\mbox{Rank} (\Delta^\omega_{G_n}-\Delta^\omega_{G'_n})<2\e|E(G_n)|\,.
\end{equation}
Indeed, let $1_p\in l^2(V(G_n))$ be the function, where $1_p(q)=0$ if $p\neq q$
and $1_p(p)=1$. Then $(\Delta^\omega_{G_n}-\Delta^\omega_{G'_n})(1_p)=0\,$
if the edges incident to $p$ are the same in $G_n$ as in $G'_n$.
Therefore
$$\dim_{\bR}\mbox{Ker}(\Delta^\omega_{G_n}-\Delta^\omega_{G'_n})\geq
|V(G_n)|-2\e|E(G_n)|\,.$$
Consequently (\ref{rank}) holds. Hence the lemma follows immediately
from Lemma 3.5 \cite{Elek2}.
Now we prove the analog of Lemma \ref{copy}. 
\begin{lemma}\label{lem2}
Let the finite graph $A$ be the disjoint union
of $m_1$ copies of $A_1$, $m_2$ copies of $A_2$,\dots, $m_l$ copies
of $A_l$ as in Lemma \ref{copy}.
Then
$$N_{\Delta^\omega_A}(\lambda)=\sum_{i=1}^l
w_i N_{\Delta^\omega_{A_i}}(\lambda)\,,$$
where $w_i=\frac{m_i |V(A_i)|}{\sum_{i=1}^l m_i |V(A_i)|}\,.$
\end{lemma}
\proof
Clearly,
$s_{\Delta^\omega_H}(\lambda)=\sum_{i=1}^l m_i 
s_{\Delta^\omega_{A_i}}(\lambda)\,.$ Therefore
$$N_{\Delta^\omega_H}(\lambda)=
\frac{\sum_{i=1}^l m_i 
s_{\Delta^\omega_{A_i}}(\lambda)}{\sum_{i=1}^l m_i |V(A_i)|}=
\sum_{i=1}^l
w_i N_{\Delta^\omega_{A_i}}(\lambda)\,.\quad\qed $$

\vskip0.2in
\noindent
For each $\omega\in\Omega_G^X$ we have a natural vertex-labeling of $G$
(and of its subgraphs), where $\omega(p)$ is the label of the vertex $p$.
Clearly, if $\omega$ and $\omega'$ coincide on the finite graph $F$
then $\Delta^\omega_F= \Delta^{\omega'}_F$.
Now let $H$ be a finite simple graph, $|V(H)|\leq K$, where $K$ is the
constant in Theorem \ref{main}. Let
$\{H_\alpha\}_{\alpha\in I_H}$ be the set
of all vertex-labellings of $H$ by $\{r_1,r_2,\dots,r_m\}$ up to
labelled-isomorphisms.
By the law of large numbers, for almost all $\omega\in\Omega^X_G$
the number of labelled vertices in $G_n'$ belonging to a vertex-labelled
component labelled-isomorphic to $H_\alpha$ divided by $|V(G_n)|$
converges to a constant $q(H_\alpha)$. Notice that
$q(H_\alpha)=c_H p(H_\alpha)$, where $p(H_\alpha)$ is the
probability that a $X$-random labelling of the vertices of $H$ is
labelled-isomorphic to $H_\alpha$.

\noindent
Hence by Lemma \ref{lem2}, for almost all $\omega\in\Omega^X_G$
the functions
$\{N_{\Delta^\omega_{G_n'}}\}^\infty_{n=1}$ converge uniformly
to a function $N_\e$ that does not depend on $\omega$. 
Let $W\subset\Omega_G^X$ be the set of elements such
that $\{N_{\Delta^\omega_{G_n'}}\}^\infty_{n=1}$ converge uniformly
to $N_{\frac{1}{k}}$ for any $k\geq 1$. Clearly, $\nu_G(W)=1$.
The following lemma finishes the proof of our Theorem.
\begin{lemma}
$\{N_{\frac{1}{k}}\}^\infty_{k=1}$ converge uniformly to a function 
$N_G^X$. Also,
for any $\omega\in W$
$\{N_{\Delta^\omega_{G_n}}\}^\infty_{n=1}$ converge uniformly to 
$N_G^X$.
\end{lemma}
\proof
By Lemma \ref{analoglemma} if $\omega\in W$ then for large enough $n$
$$|N_{\Delta^\omega_{G_n}}(\lambda)- N_{\frac{1}{k}}(\lambda)|\leq
2d\frac{1}{k}\,\, \mbox{for} \,\,-\infty<\lambda<\infty\,.$$
Therefore $\{N_{\frac{1}{k}}\}^\infty_{k=1}$ form a Cauchy-sequence and
consequently 
$\{N_{\Delta^\omega_{G_n}}\}^\infty_{n=1}$ converge uniformly to 
$N_G^X$. \qed

\subsection{Independence ratio, entropy and log-partitions}
First recall the notion of some graph parameters associated to independent
sets and matchings. Let $H$ be a finite graphs.
\begin{itemize}
\item
Let $I(H)$ be the maximal size of an independent subset in $V(H)$.
The number $\frac{I(H)}{|V(H)|}$ is called the {independence ratio} of $H$.
\item Let $M(H)$ be the maximal size of a matching in $E(H)$.
The number $\frac{M(H)}{|V(H)|}$ is called the {matching ratio} of $H$.
\item
Let $\pi^I_H(\lambda)=\sum_{\{S\subset V(H),\,S\, \mbox{is independent}\}}
\lambda^{|S|}$ be the partition function corresponding to the system of 
independent subsets.
\item
Let $\pi^M_H(\lambda)=\sum_{\{T\subset E(H),\,T\, \mbox{is a matching}\}}
\lambda^{|T|}$ be the partition function corresponding to the system of 
matchings.
\end{itemize}
Now we prove that all the graph parameters above are continuous in $\grf$.
That is we show the analog of the existence of the integrated
density of states for the quantities above.
\begin{theorem} \label{utolsotetel}
Let $G$ be an infinite connected graph such that
$|B_r(x)|\leq f(r)$, for any $x\in V(G)$ (where $f$ is of subexponential
growth)
with UPF and $\{G_n\}^\infty_{n=1}$ be a \Fo-sequence. Then
\begin{description}
\item{(a)} $\limn\frac{I(G_n)}{|V(G_n)|}$ exists.
\item{(b)} $ \limn\frac{M(G_n)}{|V(G_n)|}$ exists.
\item{(c)} $\limn \frac{\log \pi^I_{G_n}(\lambda)} {|V(G_n)|}$ exists
for all $0<\lambda<\infty\,.$
\item{(d)}$ \limn \frac{\log \pi^M_{G_n}(\lambda)} {|V(G_n)|}$ exists
for all $0<\lambda<\infty\,.$
\end{description}
Note that if $\lambda=1$ then the limit value is the associated entropy.
\end{theorem}
\proof
We prove only (a) and (c), since the proofs of (b) and (d) are completely
similar. Let $\e>0$ and $G_n'\subset G_n$ be the spanning graphs as in
Theorem \ref{main}. Again, we prove the continuity lemma.
\begin{lemma}
\label{dudo3}
$$\left|\frac{I(G'_n)}{|V(G_n)|}- \frac{I(G_n)}{|V(G_n)|}\right|< \e d $$
$$\left|\frac{\log \pi^I_{G'_n}(\lambda)} {|V(G_n)|}-
\frac{\log \pi^I_{G_n}(\lambda)} {|V(G_n)|}\right|<(
\log (\max(1,\lambda))+2)\e d\,.$$
\end{lemma}
\proof Clearly, $I(G'_n)\geq I(G_n)$.
Let $A$ be a maximal independent subset of $G_n'$. Then if we delete the
vertices from $A$ that are incident to an edge of
$E(G_n)\backslash E(G_n')$, the remaining set is an independent subset of
the graph $G_n$. Hence
$$\left|\frac{I(G'_n)}{|V(G_n)|}- \frac{I(G_n)}{|V(G_n)|}\right|< \e d \,.$$
Now let $S$ be an independent subset of the graph $G_n$. Denote by 
$Q(S)$ the set of independent subsets $S'$ of $G'_n$ such that
$S\subseteq S'$ and if $p\in S'\backslash S$ then $p$ is incident
to an edge of $E(G_n)\backslash E(G_n')\,.$
Observe that
$$
\sum_{\{S'\subset V(G_n)\,\mid\,S'\, \mbox{is independent in}\, G'_n\}}
\lambda^{|S'|} \leq 
\sum_{\{S\subset V(G_n)\,\mid\,S\,  \mbox{is independent in}\, G_n\}} 
\sum_{S'\in Q(S)} \lambda^{|S'|}
\leq$$
$$\leq \sum_{\{S\subset V(G_n)\,\mid\,S\, 
 \mbox{is independent in}\,\, G_n\}} \lambda^{|S|}
\max(1,\lambda)^{\e d |V(G)|} 2^{\e d |V(G)|}\,.$$
That is
$$\log \pi^I_{G'_n}(\lambda)\leq \log \pi^I_{G_n}(\lambda)+
(\log(\max(1,\lambda))+2)\e d |V(G_n)|\,.\quad\qed$$
\begin{lemma}
\label{dudo4}
Let $A_1, A_2, A_3,\dots, A_l$ be finite simple graphs. Suppose that
the graph $H$ consists of $m_1$ disjoint copies of $A_1$ and $m_2$
disjoint copies of $A_2$, \dots and $m_l$ disjoint copies of $A_l$.
Then
\begin{equation} \label{eb1}
\frac{I(H)}{|V(H)|}=\sum^l_{i=1} w_i \frac{I(A_i)}{|V(A_i)|}
\end{equation}
\begin{equation} \label{eb2}
\frac{\log \pi^I_H (\lambda)}{|V(H)|}=\sum^l_{i=1} w_i
\frac{\log \pi^I_{A_i} (\lambda)}{|V(A_i)|}\,,
\end{equation}
where $w_i=\frac{m_i |V(A_i)|}{\sum_{i=1}^k m_i |V(A_i)|}\,.$
\end{lemma}
\proof
Note that $I(H)=\sum^l_{i=1} m_i I(A_i)$ and
$\pi^I_H(\lambda)=\prod_{i=1}^l (\pi^I_{A_i})^{m_i}$. That is
$\log(\pi^I_H(\lambda))=\sum_{i=1}^l m_i \log (\pi^I_{A_i})\,.$
Now we proceed as in Lemma \ref{analoglemma}\qed.

\noindent
By the two preceding lemmas Theorem \ref{utolsotetel} easily follows. \qed

\vskip 0.2in
\noindent
{\bf Remark:} In Theorem 3. \cite{BG} the authors proved that 
$\limn\frac{I(G_n)}{|V(G_n)|}$ exists if $\{G_n\}^\infty_{n=1}$ 
is a $r$-regular large girth sequence, where $2\leq r\leq 5$.

\end{document}